\documentclass[12pt]{amsart}
\advance\hoffset by -0.5 truecm
\usepackage{amssymb}
\newtheorem{Theorem}{Theorem}[section]
\newtheorem{Lemma}[Theorem]{Lemma}
\newtheorem{Corollary}[Theorem]{Corollary}

\newtheorem{Proposition}[Theorem]{Proposition}

\newtheorem{Definition}[Theorem]{Definition}
\newtheorem{Remark}[Theorem]{Remark}

\def \dim{{\mbox {dim}}\,}

\def\V{\mbox{Var}}

\def\Z{{\mathbb Z}}
\def\R\re
\def\V{\bf V}

\def \la{\lambda}

\def \re{{\mathbb R}}
\def \Q{{\mathbb Q}}

\def \T{{\mathbb T}}
\def \C{{\mathbb C}}
\def \M{{\widetilde{M}}}

\def \H{{\mathbb H}}

\def \0{\lambda_{0}}
\def \la{\lambda}
\def \ga{\gamma}

\def\h{{\rm h}_{\rm top}(g)}

\def \Eu{{\mathbb E}}

\begin{document}
\title[Entropy and collapsing of complex
surfaces]{Entropy and collapsing of compact
complex surfaces}

\author[G. P. Paternain]{Gabriel P. Paternain}\thanks{G. P. Paternain
 was partially
 supported by CIMAT, Guanajuato, M\'exico}
 \address{ Department of Pure Mathematics and Mathematical Statistics,
University of Cambridge,
Cambridge CB3 0WB, England}
 \email {g.p.paternain@dpmms.cam.ac.uk}

\author[J. Petean]{Jimmy Petean}
 \address{CIMAT  \\
          A.P. 402, 36000 \\
          Guanajuato. Gto. \\
          M\'exico.}
\email{jimmy@cimat.mx}

\thanks{J. Petean is supported by grant 37558-E of CONACYT}


\date{revised and improved version: June 2003}


\begin{abstract} We study the problem of existence
of $\mathcal{F}$-structures (in the sense of Cheeger and Gromov,
but not necessarily polarized) on compact
complex surfaces. We give a complete classification
of compact complex surfaces of K\"{a}hler type
admitting  $\mathcal{F}$-structures. In the
non-K\"{a}hler case we give a complete classification
modulo the gap in the classification of
surfaces of class VII. In all these examples a surface admits an
$\mathcal{F}$-structure if and only if it admits a $\mathcal{T}$-structure.

We then use these results to study the minimal entropy
problem for compact complex surfaces: we prove that, modulo
the gap in the classification of surfaces of class VII, all
compact complex surfaces of Kodaira dimension $\leq 1$ have
minimal entropy 0 and such a surface admits a smooth metric $g$
with $\h=0$  if and only if it is
$\C P^2$, a ruled surface of genus $0$ or $1$,
a Hopf surface, a complex torus, a Kodaira surface,
a hyperelliptic surface or a Kodaira surface modulo
a finite group. The key result we use to prove this, is 
a new topological obstruction to the existence of metrics
with vanishing topological entropy.

Finally we show that these results fit perfectly into
Wall's study of geometric structures on compact complex surfaces.
For instance, we show that the minimal entropy problem 
can be solved for
a minimal compact K\"{a}hler surface $S$ of Kodaira dimension
$-\infty$, 0 or 1 if and only if $S$ admits a geometric structure
modelled on $\C P^2$, $S^2\times S^2$, $S^2\times\Eu^2$ or 
$\Eu^4$. 
\end{abstract}

\maketitle

\section{Introduction}
The concept of an  $\mathcal{F}$-structure was introduced
by M. Gromov in \cite{Gromov} as a natural
generalization of a torus action on a manifold. Then
J. Cheeger and Gromov 
(\cite{Gromov, CG1, CG2}) related the existence of
 $\mathcal{F}$-structures with the possibility
of {\it collapsing} a manifold with bounds on certain
geometric invariants (curvature, injectivity radius,
diameter, etc.). These results require the
 $\mathcal{F}$-structures to be of a very restrictive
type, namely, {\it polarized} $\mathcal{F}$-structures.
 An $\mathcal{F}$-structure is a sheaf of tori
acting on open subsets of a manifold (with certain
compatibility on the intersection of the open subsets);
polarized means (broadly speaking) that the
dimension of the orbits of
the actions is locally constant. In \cite{PP} the authors
proved that if a closed, smooth manifold $M$ admits
a general  $\mathcal{F}$-structure, then one can collapse
the volume of $M$ with curvature bounded from below and
with bounded entropy. Namely, there exists a sequence
$g_i$ of Riemannian metrics on $M$ such that there is
a uniform upper bound on the entropy of the geodesic 
flow of $g_i$
and a uniform lower bound for the sectional curvatures
of $g_i$, while the volume of $(M,g_i )$ converges to 0.
This is our main motivation for looking at 
the problem of existence of $\mathcal{F}$-structures
on a given manifold.

On the other hand, compact complex surfaces provide
the most important family of examples of smooth closed
4-manifolds and they are fundamental in the study of
differential topology in dimension 4. Therefore
it seems worthwhile to investigate which of them admit  general
$\mathcal{F}$-structures. There is a well-known
classification of compact complex surfaces
(the Enriques-Kodaira classification). The
classification first divides them
according to the Kodaira dimension $\kappa$, which could be
$-\infty ,0,1$ or $2$. We briefly recall its definition.
Let $S$ be a compact complex surface and
let $P_m(S)$ be the dimension of the space of holomorphic sections
of the $m$th tensor power of the canonical line bundle of $S$.
If $P_m(S)=0$ for all $m\geq 1$, then $\kappa=-\infty$, otherwise,
$\kappa=\limsup_{m\to\infty}\frac{\log P_m(S)}{\log m}$.

Surfaces of Kodaira dimension
0 or 1 are all diffeomorphic to elliptic
surfaces. Surfaces of Kodaira dimension 2
are known as {\it surfaces of general type} since
these are the {\it generic} surfaces. It is for
surfaces of Kodaira dimension $-\infty$ that the
classification is not complete (at least not
known to be complete). K\"{a}hler surfaces of
Kodaira dimension $-\infty$ are obtained by blowing up
$\C P^2$ or $\C P^1$-bundles over a Riemann
surface and are perfectly understood. The gap in
the classification appears in the non-K\"{a}hler case.
Non-K\"{a}hler surfaces of Kodaira dimension $-\infty$
are called surfaces of class VII. There are quite a few examples.
First, we have the Hopf surfaces, characterized
by the fact that their universal cover is $\C^2
 - \{ 0 \}$. Then
we have the Inoue surfaces with vanishing second Betti
number and some compact elliptic surfaces.
Finally we have the surfaces with {\it global
spherical shells}, which have positive second Betti
number. These are all the known examples and it has
been conjectured
that these are actually all the
minimal surfaces of class VII. See Section 2 for more details.
The first aim of this article is to prove the following
results:

\medskip

\noindent {\bf Theorem A.} {\it Let $M$ be a compact K\"{a}hler
surface. Then $M$ admits an
$\mathcal{F}$-structure if and
only if the Kodaira dimension of $M$ is different
from $2$. Actually, K\"{a}hler surfaces of Kodaira
number $-\infty , 0 $ or $1$ admit $\mathcal{T}$-structures.}

\medskip

\noindent {\bf Theorem B.} {\it All known examples of compact complex surfaces
which are not of K\"{a}hler type admit
$\mathcal{T}$-structures.}

\medskip
 
Using the result of \cite{PP} mentioned above we obtain:

\medskip

\noindent {\bf Corollary.} {\it Every compact complex surface which is
not of general type has minimal entropy 0
and collapses with sectional curvature
bounded from below, except, perhaps,
some new examples of surfaces of class VII.}

\medskip

We will prove these results in Sections 2 and 3.
Next, we turn our attention to the problem
of determining for which compact complex surfaces the
minimal entropy problem can be solved. Recall that the
minimal entropy of a closed manifold $M$, ${\bf h} (M)$,
is the infimum of the topological entropy of $C^{\infty}$ Riemannian
metrics on $M$ with volume one. We say that the minimal
entropy problem for $M$ can be solved if there exists a $C^{\infty}$
metric $g$ with volume one for which $\h ={\bf h} (M)$.
The guiding principle is that manifolds for which the
minimal entropy problem can be solved should be topologically simple
and/or support geometric structures. See \cite{AP,PP} for details and related
results.
In Section 5 we begin the study of the minimal entropy problem
for compact complex surfaces of Kodaira dimension 
different from 2. From the Corollary after Theorem B it follows
that ${\bf h} (M)=0$ for all such surfaces, modulo
the gap in the classification of surfaces of class VII. Therefore
the minimal entropy problem can be solved for these surfaces
only if they admit a metric with vanishing topological entropy.

Finding topological obstructions to the existence of metrics with
zero entropy is a subtle problem. Up to now, the known results 
were of two kinds. Either the fundamental group was big in the 
sense that it had exponential growth, or the manifold was simply 
connected and the exponential
growth in the topology was located in the loop space homology.
An important ingredient in the second case is a theorem of Gromov 
\cite{G1, bates} which asserts that on a closed simply connected
manifold $M$, homology classes of the pointed loop space $\Omega M$
can be represented by cycles formed with paths of appropriately
bounded length. In the appendix we will show that Gromov's 
theorem also works in the case that $M$ is compact with
non-empty boundary. This will help us to deal with the case of
manifolds with an  infinite fundamental group of subexponential
growth.  Essentially nothing was known in this case. 
Many compact complex surfaces fall into this category 
(as well as many other manifolds!). 
More precisely, we will prove the following result. (As usual, given a space
$Y$ and points $x,y\in Y$, $\Omega(Y,x,y)$ will denote the space of paths from $x$ to $y$.)

\medskip

\noindent {\bf Theorem C.} {\it Let $M$ be a closed manifold and let $\M$ be its universal covering.
Let $X\subset\M$ be a compact simply connected submanifold, 
possibly with boundary. Consider
points $x,z \in X$ and $y\in \M$ and a path $\alpha$ from $z$ to $y$.
Define a map $\iota: \Omega (X,x,z)\to \Omega (\M,x,y)$
by $\iota(\tau)=\tau *\alpha$ and suppose $\iota$
induces a monomorphism in homology with coefficients in the prime field $k_p$, $p$ prime or zero. Then, for any $C^{\infty}$ 
Riemannian metric $g$ on $M$ there exists a positive constant $C(g)$ such that
\[\h\geq \frac{\la(g)}{2}+C(g)\,(-\log R_{\Omega X,p}),\]
where $ R_{\Omega X,p}$ is the radius of convergence of the
Poincar\'e series
$$\sum_{i\geq 0}b_{i}(\Omega X,k_p)t^{i}.$$
}

\medskip

The quantity $\la(g)$ is the {\it volume entropy} of the manifold and it
is defined as the exponential growth rate of the volume of balls in $\M$.
Manning's inequality
\cite{Ma} asserts that for any metric $g$,
${\rm h}_{\rm top}(g)\geq \la(g)$.
It is well known that $\la(g)>0$ if and only if $\pi_1(M)$ has exponential
growth. Theorem C does not really say much if $\la(g)>0$, but it is most
interesting when $\la(g)=0$ and $R_{\Omega X,p}<1$ for some $p$.

We describe a noteworthy application of Theorem C.

\medskip

\noindent {\bf Theorem D.} {\it Let $M$ be a closed manifold of dimension $n\geq 3$.
Suppose that $M$ can be decomposed as $X_{1}\# X_{2}$, where 
the order of the fundamental group of $X_1$ is at least $3$.
If $M$ admits a $C^{\infty}$ Riemannian metric with zero
topological entropy, then $X_2$ is a homotopy sphere.}

\medskip

Theorems C and D will be proved in Section 4. Note that Theorem D is
optimal. B. Totaro proved in \cite{T} that $\re P^{n}\#\re P^{n}$
and $\C P^{n}\# \re P^{2n}$ are diffeomorphic to biquotients and hence
they have real analytic metrics with zero topological entropy.

We can now state the theorems that will give 
complete solutions for the minimal
entropy problem for compact complex surfaces of Kodaira dimension 
different from 2. We phrase the solution of the minimal entropy 
problem in terms 
of 4-dimensional geometric structures. We say that $M$ admits a
geometric structure if $M$ admits a locally homogeneous metric. 
This means that if we endow $\M$ with the corresponding covering 
metric, then $\M$ is a complete homogeneous 
space and we say that $M$ admits a geometric structure modelled 
on $\M$.
In dimension 4, the maximal geometric models have been classified by 
R.O. Filipkiewicz \cite{Fi}. C.T.C Wall \cite{wall} studied the relationship 
between geometric structures and complex structures for compact complex 
surfaces. Wall's results will be quite useful for the proof of the
two theorems below.

\medskip

\noindent {\bf Theorem E.} {\it Let $S$ be a compact complex surface not of K\"ahler type.
 Modulo the gap in the classification of class VII surfaces we have:
The minimal entropy of $S$ is zero and the following are equivalent:
\begin{enumerate}
\item The minimal entropy problem can be solved for $S$;
\item $S$ admits a smooth metric $g$ with $\h=0$;
\item $S$ admits a geometric structure modelled on $S^{3}\times \Eu^1$ 
or {\rm Nil}$^3\times \Eu^1$;
\item $S$ has $\kappa=-\infty,0$ and is diffeomorphic to one of the following:
a Hopf surface, a Kodaira surface, or a Kodaira surface modulo a finite group.
\end{enumerate}
}

\medskip

Recall that a Kodaira surface is a surface $S$ with $K_{S}=\mathcal O_{S}$
and $b_1=3$. Such a surface is diffeomorphic to the product $S^{1}\times N$, where $N$
is a 3-dimensional nilmanifold.

In the K\"ahler case we have:

\medskip

\noindent {\bf Theorem F.} {\it Let $S$ be a compact complex K\"ahler surface with Kodaira dimension
$\kappa\leq 1$. We have:
The minimal entropy of $S$ is zero and the following are equivalent:
\begin{enumerate}
\item The minimal entropy problem can be solved for $S$;
\item $S$ admits a smooth metric $g$ with $\h=0$;
\item $S$ admits a geometric structure modelled on $\C P^2$, $S^2\times S^2$, $S^2\times\Eu^2$
 or $\Eu^4$ or $S$ is diffeomorphic to $\C P^2\#\overline{\C P^{2}}$;
\item $S$ has $\kappa=-\infty,0$ and is diffeomorphic to one of the following:
$\C P^2$, a ruled surface of genus $0$ or $1$, a complex torus or a hyperelliptic surface.
\end{enumerate}
}

Note that $\C P^2\#\overline{\C P^{2}}$ is the only non-minimal surface that admits
a metric of zero entropy. It is also the only surface that admits a metric with
zero entropy and no geometric structure.

\medskip

{\it Acknowledgements:} We thank B. Totaro for useful discussions. We
also thank C. LeBrun for useful comments and a very helpful hint
in the construction of retractions in Section 4.
The first author thanks the CIMAT, Guanajuato, M\'exico for hospitality
and support while part of this work was carried out.

\section{${\mathcal T}$-structures and compact complex surfaces}

We will first review the definition of a ${\mathcal T}$-structure
and the results we will need about them,
and then study the problem of existence of such
structures on compact complex surfaces using the
Enriques-Kodaira classification.
We will prove here Theorems A and B 
in the introduction modulo
the construction of ${\mathcal T}$-structures on surfaces of
Kodaira dimension $-\infty$ that we will carry out in the
next section.

\begin{Definition} A ${\mathcal T}$-structure on a smooth
closed manifold $M$ is a finite
open cover $(U_i )_{i=1,...,l}$ of $M$ and
a non-trivial torus action on each $U_i$ such that the
intersections of the open subsets are invariant
(through all the corresponding torus actions) and
the actions commute.
\end{Definition}

The ${\mathcal T}$-structure is called {\it polarized}
if the torus actions on each $U_i$
are locally free and in the
intersections the dimension of the orbits (of
the corresponding torus action) is constant.
The structure is called {\it pure} if the dimension
of the orbits is constant.

\begin{Remark} {\rm If two manifolds of dimension greater than 2
admit ${\mathcal T}$-structures then their
connected sum also admits one (see \cite{Soma,
Gromov} for the case of polarized structures on
odd dimensional manifolds and \cite{PP} for
the general case). We will use this result mainly
to reduce the problem of existence of a ${\mathcal T}$-structure
on a compact complex surface $M$ to the case when
$M$ is minimal, since non-minimal surfaces are
obtained from minimal ones by taking connected sums
with $\overline{{\C P}^2}$.}
\end{Remark}

We will also use the following two results from
\cite{PP}:

\begin{Theorem} If a closed smooth manifold
$M$ admits a ${\mathcal T}$-structure then the minimal
entropy of $M$ is 0 and $M$ collapses with curvature
bounded below.
\end{Theorem}

\noindent
(Actually, the conclusion of the theorem still holds
if one only assumes the existence of an
$\mathcal{F}$-structure.)

\begin{Theorem} Every elliptic compact complex
surface admits a ${\mathcal T}$-structure.
\end{Theorem}

In the next section we will show that Hopf surfaces,
Inoue surfaces with $b_2$=0 and ruled surfaces all
admit ${\mathcal T}$-structures. Using this, we are in a position
to prove Theorems A and B.

Minimal compact complex surfaces of K\"{a}hler type
are: ${\C P}^2$, ruled surfaces, surfaces which are 
deformation equivalent to 
elliptic surfaces
(with even first Betti number) and surfaces of general
type. From Theorem 2.3 we know that if a closed
manifold admits a ${\mathcal T}$-structure then it collapses
with curvature bounded below, and therefore its Yamabe
invariant is non-negative (see \cite{PP}). But it is known that
surfaces of general type have strictly negative Yamabe
invariant (see \cite{LeBrun}, for instance). 
Therefore they cannot admit ${\mathcal T}$-structures.
All the other cases are covered by the previous
comments and Theorem A is proved.

Minimal compact complex surfaces which are not
of K\"{a}hler type are either elliptic or
surfaces of class VII. We therefore only need to consider
surfaces of class VII. These are the surfaces which have  first
Betti number equal to 1 and  Kodaira dimension equal to
$-\infty$. They were first studied by Kodaira
\cite{Kodaira}. All the known examples are the following:

i) Elliptic surfaces of class VII;

ii) Hopf surfaces;

iii) Inoue surfaces with $b_2 =0$;

iv) Surfaces with a global spherical shell.

Elliptic surfaces admit ${\mathcal T}$-structures by Theorem 2.4.
We will exhibit ${\mathcal T}$-structures on Hopf surfaces and Inoue
surfaces with $b_2 =0$ in the following section.
Surfaces with a global spherical
shell have positive second Betti number and
are all diffeomorphic to a connected sum of
$S^1 \times S^3$ with $b_2 (M)$ copies of
$\overline{{\C P}^2}$. Therefore they also admit 
${\mathcal T}$-structures by Remark 2.2.

Modulo the constructions in the following section, this
finishes the proof of Theorem B.

\begin{Remark}{\rm  Kodaira proved that a surface of class VII
with vanishing second Betti number which contains at least a
curve is either an elliptic surface or a Hopf surface.
Later on Inoue \cite{Inoue} constructed his examples
of surfaces of class VII with $b_2 =0$ which contain no
curves. He also proves that if a surface of class VII
with $b_2 =0$ and no curves admits a line bundle
$F_0$ such that $\mbox{\rm dim}\, H^0 ({\Omega}^1 (F_0 ))>0$
then the surface is one of his examples. Later
F.A. Bogomolov \cite{Bogomolov1, Bogomolov2}
gave a proof that such a
line bundle always exists, but apparently his
proof was difficult to understand.
Later, J. Li, S. T. Yau and F. Zheng \cite{LYZ1, LYZ2}
gave a much simpler proof.
This result concluded
the classification of  surfaces of class VII with $b_2 =0$.
Surfaces with a global spherical shell have positive
second Betti number and it has been conjectured
\cite[page 220]{Nakamura} that every surface of class VII
with positive second Betti number admits a global
spherical shell. If this conjecture were true, the classification
would be complete and Theorem B would cover all non-K\"{a}hler
surfaces.}
\end{Remark}

\section{Construction of ${\mathcal T}$-structures}

In this section we will carry out the construction
of ${\mathcal T}$-structures on Hopf surfaces, Inoue surfaces with vanishing second Betti
number and ruled surfaces.

\subsection{Hopf surfaces}
The topology of Hopf surfaces
was studied by M. Kato  \cite[Theorem 9]{Kato}.
There he proves that every Hopf surface is diffeomorphic
either to a product $S^1 \times (S^3 /H)$ or to an
$S^3 /H$-bundle over $S^1$ with an involution as
the transition function (here $H$ is some subgroup of
$GL(2,{\C})$).
In both cases the surface has
an obvious locally free
(or free) $S^1$-action.

Therefore we can write

\begin{Lemma} Every Hopf surface admits a locally free
$S^1$-action.
\end{Lemma}

\subsection{Inoue surfaces with $b_2 =0$}
There are three classes of Inoue surfaces with
vanishing second Betti number; named $S_M$,
$S_{N,p,q,r;t}^+$, and $S_{N,p,q,r}^-$. They were discovered
by Inoue in \cite{Inoue}.

\vspace{.5cm}

We first consider the surfaces $S_M$.
Here $M$ denotes a $3\times 3$
integral matrix
with determinant 1 and eigenvalues $\alpha$, $\beta$
and $\overline{\beta}$ with $\alpha >1$ and $\beta
\neq \overline{\beta}$. $M$ induces a transformation
$f_M$ of the torus $T^3 = {\re}^3 / {\Z}^3$ (sending the class of
$x\in {\re}^3 $ to the class
of $Mx \in {\re}^3$) and the surface
$S_M$ is diffeomorphic to the $T^3$-bundle over
$S^1$ with $f_M$ as transition function. Then $S_M$
admits a pure polarized ${\mathcal T}$-structure of rank 3, whose
orbits are the fibers of the fibration.

\vspace{.5cm}

Now consider the surfaces $S_{N,p,q,r;t}^+$.
Here $N=(n_{ij}) \in SL(2,{\Z})$ is a
unimodular integral $2\times 2$
matrix with  real eigenvalues $\alpha$ and
${\alpha}^{-1}$, with $\alpha >1$. Then pick any integers
$p,q$ and $r$ with $r\neq 0$ and a complex number $t$.
These are the parameters for the surfaces
$S_{N,p,q,r;t}^+$. Then one picks eigenvectors
$(a_1 ,a_2 )$ and $(b_1 ,b_2 )$ (with eigenvalues $\alpha$
and ${\alpha}^{-1}$, respectively) and considers
the analytic automorphisms of ${\mathbb H}\times {\C}$
given by

$$ g_0 (w,z)=(\alpha w,z+t);$$

$$ g_i (w,z)= (w+a_i , z+b_i w +c_i ), \ \ \ i=1,2; $$

$$ g_3 (w,z)=\left( w, z+(b_1 a_2 -b_2 a_1 )/r \right), $$

\noindent
where $c_1 $ and $c_2$ are constants which are obtained
from $N$ and its eigenvectors 
(see \cite{Inoue} for the formulas)
and ${\mathbb H}=\{ w\in
{\C} : \mbox{\rm Im}(w) > 0\}$. The surface $S_{N,p,q,r;t}^+$ is
then the quotient of ${\mathbb H}\times {\C}$ by the group
of analytic automorphisms generated by $g_0 ,g_1 ,g_2$ and
$g_3$.

Think of $S^1$ as obtained from the interval $[0,(b_1 a_2
-b_2 a_1 )/r ]$ by identifying its endpoints
(of course, $(1/r)(b_1 a_2 -b_2 a_1)$ can be assumed
to be positive). Then it is
easy to check that the formula

$$ (\theta ,w,z) \mapsto (w,\theta +z)$$

\noindent
defines a locally free $S^1$-action on $S_{N,p,q,r;t}^+$.

\vspace{.5cm}

Finally, let us consider the surfaces $S_{N,p,q,r}^-$. Here
$N\in GL(2,{\Z})$ is an integral $2\times 2$ matrix with
determinant equal to $-1$ having real eigenvalues
$\alpha >1$ and $-{\alpha}^{-1}$. Choose non-trivial
eigenvectors $(a_1 ,a_2 )$ and $(b_1 ,b_2 )$ corresponding
to $\alpha$ and  $-{\alpha}^{-1}$ respectively. Fix integers
$p,q$ and $r$, with $r\neq 0$, and consider the
analytic automorphisms of ${\mathbb H} \times {\C}$ given
by

$$g_0 (w,z)= (\alpha w,-z);$$

$$g_i (w,z)= (w+a_i ,z+b_i w + c_i ), \ \ \ i=1,2;$$

$$g_3 (w,z)= \left( w,z+(b_1 a_2 -b_2 a_1 )/r \right) .$$

\noindent
Here $c_1$ and $c_2$ are constants obtained from
$N, a_1 ,a_2 ,b_1 ,b_2 ,p$ and $q$ (see \cite{Inoue}
for the formulas). Then
$S_{N,p,q,r}^-$ is the quotient of ${\mathbb H}\times {\C}$
by the group of automorphisms generated by the $g_i$'s.

Consider the map
${\mathbb H}\times {\C} \rightarrow {\re}_{>0}$ given
by the projection onto the second coordinate (the
imaginary part of the first complex coordinate).
It is easy to check that it induces a map
$S_{N,p,q,r}^-  \rightarrow S^1$, where $S^1$ is the
interval $[1,\alpha ]$ with its endpoints identified.

Consider the group  $G_y$ of analytic automorphisms
generated by $g_1 ,g_2$ and $g_3$. The projection is invariant
under these automorphisms and the fibers of
it are identified with the quotient of ${\re}^3$
by the action of this group (with $y=\mbox{\rm Im}(w)$ fixed).

There is an $S^1$-action on each fiber given by
$\theta . \ (\mbox{\rm Re}(w),z)=(\mbox{\rm Re}(w),\mbox{\rm Re}(z)+\theta ,\mbox{\rm Im}(z))$,
where $S^1$ is now thought of as the interval
 $[0,(b_1 a_2
-b_2 a_1 )/r ]$ with its endpoints identified.

Now $S_{N,p,q,r}^-  \rightarrow S^1$ gives  $S_{N,p,q,r}^-$
a fiber bundle structure. On each fiber we have
the $S^1$-action mentioned above.
Call $._1$ this action. Under the transition map this
action translate to a new action which we will call
$._2$. If these two actions were to coincide we would have a
well defined $S^1$-action, as for the surfaces
$S_{N,p,q,r;t}^+$. But in order to have a
${\mathcal T}$-structure we only need these actions to
commute. Now a direct computation shows that

$$ \theta ._2 (\mbox{\rm Re}(w),\mbox{\rm Re}(z),\mbox{\rm Im}(z)) =
 (-\theta )._1  (\mbox{\rm Re}(w),\mbox{\rm Re}(z),\mbox{\rm Im}(z)) .$$

\noindent
Therefore  $S_{N,p,q,r}^- $ admits a pure polarized
${\mathcal T}$-structure of rank one.

Summarizing we can write:

\begin{Lemma} Every Inoue surface with vanishing second
Betti number admits a pure polarized
${\mathcal T}$-structure of positive rank.
\end{Lemma}

\subsection{Ruled surfaces} Every minimal
ruled surface $M$ is
diffeomorphic to the projectivization of a decomposable
vector bundle of rank 2 over a compact Riemann surface $S$.
Moreover we can assume that
one of the line bundles in the decomposition
is trivial and so there exists a complex
line bundle $L$ over $S$ such that $M$ is diffeomorphic
to ${\bf P}(L \oplus {\C})$. The $S^1$-action on
${\C}$ then induces an $S^1$-action on $M$. Therefore
we have:

\begin{Lemma}Every ruled surface admits a non-trivial
$S^1$-action.
\end{Lemma}

\section{Obstructions to zero entropy}
\label{mecanismos}

In this section we prove Theorems C and D.
We will need the following version of a result due to 
M. Gromov \cite{G1} for manifolds with non-empty boundary.
Gromov's original proof in \cite{G1} is very short. A more
detailed proof is given in \cite[p. 102]{pansu} but it seems
unclear at certain points. There is a corrected proof in the
revised an updated edition \cite{bates}. Also I. Babenko sketches
a proof in \cite{Ba} using Adams's cobar construction.
Finally there is also a proof in \cite{P2}. 
In all these references the manifold is assumed to have
empty boundary. The proofs in the cited references
work also in the case of compact manifolds with non-empty boundary
with some modifications. For completeness
we will sketch how to modify the proof given in \cite{P2}
in the appendix.

\begin{Theorem}Given a metric $g$ on a simply connected
compact manifold $X$ (possibly with boundary), there exists a
constant $C>0$ such that given any pair of points $x$
and $y$ in $X$ and any positive integer $i$, any element
in $H_{i}(\Omega(X,x,y))$ can be represented by a cycle
whose image lies in $\Omega^{Ci}(X,x,y)$. \label{coolfact}
\end{Theorem}

Let $M$ be a complete Riemannian
manifold. Given $x$ and $y$ in $M$ and $T>0$ let $n_{T}(x,y)$ be the number of geodesic segments
of length $\leq T$ connecting $x$ to $y$.
Let $B(x,T)$ be the closed ball with center at $x$ and radius
$T$. Clearly $n_{T}(x,y)\geq 1$ if and only if $y\in B(x,T)$. Moreover it is easy to see that
\[\int_{M}n_{T}(x,y)\,dy=\int_{B(x,T)}n_{T}(x,y)\,dy=\int_{B(0,T)}|\det d_{\theta}\exp_{x}|\,d\theta.\]

 Let $(M,g)$ be a closed Riemannian manifold and let
$\widetilde{M}$ be its universal covering endowed with the induced
metric.  Given $x\in\widetilde{M}$, let $V(x,r)$ be the volume of
the ball with center $x$ and radius $r$.  Set
\[\lambda(g):=\lim_{r\to +\infty}\frac{1}{r}\log\,V(x,r).\]
Manning \cite{Ma} showed that the limit exists and it is independent
of $x$.


It is well known \cite{Mil} that $\lambda(g)$ is positive if and
only if $\pi_{1}(M)$ has exponential growth. Manning's inequality
\cite{Ma} asserts that for any metric $g$,
\begin{equation}
{\rm h}_{\rm top}(g)\geq \lambda(g).
\label{manineq}
\end{equation}
In particular, it follows that if $\pi_{1}(M)$ has exponential
growth then ${\rm h}_{\rm top}(g)$ is positive for any metric $g$.
This fact was first observed by E.I. Dinaburg in \cite{D}.

Given a manifold $\M$, a subset $X\subset \M$, points $x,z\in X$,
$y\in \M$ and a continuous path $\alpha$ from $z$ to $y$, there
is an inclusion $\iota:\Omega (X,x,z) \rightarrow \Omega (\M ,x,y)$ given
by $\iota(\tau)=\tau *\alpha$.
Let $I:H_*  (\Omega (X,x,z),k_p) \rightarrow 
H_* (\Omega (\M ,x,y),k_p)$ be the map induced in homology for some coefficient field $k_p$. Picking a
different path from $z$ to $y$ gives a new inclusion, homotopic to the
previous one. Therefore the map $I$ induced in homology is the same.
Moreover, it is easy to see that the question of whether the map $I$
is injective or not is independent of the choices of $x,z,y,\alpha$.
We now show:

\medskip

\noindent {\bf Theorem C.} {\it Let $M$ be a closed manifold and let $\M$ be its universal covering.
Suppose that there exists a compact simply connected submanifold
$X\subset\M$, possibly with boundary, such that for some points
$x,z\in X$, $y\in \M$ and a path from $z$ to $y$
the inclusion $\iota:\Omega (X,x,z)
\rightarrow \Omega (\M ,x,y)$ induces a monomorphism 
$I :H_* (\Omega (X,x,z),k_p)
\rightarrow H_* (\Omega (\M ,x,y),k_p)$. 
Then, for any $C^{\infty}$ Riemannian metric
$g$ on $M$ there exists a positive constant $C(g)$ such that
\[\h\geq \frac{\la(g)}{2}+C(g)\,(-\log R_{\Omega X,p}),\]
where $ R_{\Omega X,p}$ is the radius of convergence of the
Poincar\'e series
$$\sum_{i\geq 0}b_{i}(\Omega X,k_p)t^{i}.$$
}

\medskip

\begin{proof}It is well known that for any $x\in M$ (cf. \cite{P1})
\[\h\geq
\limsup_{T\to\infty}\frac{1}{T}\log\int_{M}n_{T}(x,y)\,dy.\] Let
$p:\M\to M$ be the covering projection. One easily checks that
given any $x\in \M$ we have
\[\int_{M}n_{T}(p(x),y)\,dy=\int_{\M}n_{T}(x,y)\,dy=\int_{B(x,T)}n_{T}(x,y)\,dy.\]
Thus for any $x\in \M$ we have
\begin{equation}
\h\geq
\limsup_{T\to\infty}\frac{1}{T}\log\int_{B(x,T)}n_{T}(x,y)\,dy.
\label{equis} \end{equation}

We will use the following lemma whose proof will be given after 
completing the proof of the theorem. The notation is the same as in
the theorem and the comments before it. In what follows we take $x\in X$.

\begin{Lemma} There exists a constant
$C(g)>0$ such that for any $T\geq Ci$ and any $y\in B(x,T/2)$,
every element in the image of the map $I$ 
constructed using $z=x$ and $\alpha$ a minimizing geodesic
from $x$ to $y$ can be represented by
a cycle in $\Omega^{T}(\M,x,y)$. 
\label{retraccion}
\end{Lemma}

If $x$ and $y$ are not conjugate points, then the Morse 
inequalities imply:
\[n_{T}(x,y)\geq \sum_{i\geq 0}b_{i}(\Omega^{T}(\M,x,y),k_p).\]

Using the lemma and the hypothesis that $I$ is a monomorphism,
we obtain for any positive integer $k$ and any
$y\in B(x,Ck/2)$ not conjugate to $x$,
\[n_{Ck}(x,y)\geq \sum_{i\geq 0}b_{i}(\Omega^{Ck}(\M,x,y),k_p)
\geq \sum_{i\leq k}b_{i}(\Omega(X,x,x),k_p).\] Let us now
integrate this inequality with respect to $y\in B(x,Ck/2)$.
Since $b_{i}(\Omega(X,x,x),k_p)$ is independent of $x$ we
obtain:
\[\int_{B(x,Ck/2)}n_{Ck}(x,y)\,dy\geq
V(x,Ck/2)\left(\sum_{i\leq k}b_{i}(\Omega X,k_p)\right)\] and
hence
\[C\limsup_{T\to\infty}\frac{1}{T}\log\int_{B(x,T)}n_{T}(x,y)\,dy
\geq C\frac{\la(g)}{2}+\limsup_{k\to +\infty}\frac{1}{k}\log
\sum_{i\leq k}b_{i}(\Omega X,k_p).\] 
Clearly
\[-\log\,R_{\Omega X,p}\leq
\limsup_{k\rightarrow +\infty}\frac{1}{k}\log \sum_{i=0}^{k}b_{i}(\Omega X,k_p),\]
with equality if and only if the Poincar\'e series of $\Omega X$
has infinitely many non-zero coefficientes.
Thus
\[\limsup_{T\to\infty}\frac{1}{T}\log\int_{B(x,T)}n_{T}(x,y)\,dy
\geq \frac{\la(g)}{2}+\frac{1}{C}(-\log\,R_{\Omega X,p})\] which
together with inequality (\ref{equis}) completes the proof of the
theorem.

\end{proof}

\noindent {\it Proof of Lemma \ref{retraccion}.} Consider $c\in
H_{i}(\Omega(X,x,x))$. By Theorem \ref{coolfact} there
exists a constant $C'(g)>0$ such that for any $T\geq C'i$ there
exists a cycle $\eta$ in $\Omega^{T}(X,x,x)$ such that $\eta$
represents $c$.

Let $\alpha$ be a minimizing geodesic connecting $x$ to $y$.
As before, let
$$\iota :\Omega(X,x,x)\to \Omega(\M,x,y)$$
be the map
\[\tau\mapsto \tau*\alpha.\]
Clearly the length $\ell(\iota(\tau))$ of $\iota(\tau)$ is given by
$\ell(\tau)+d(y,x)$. Let $C=2C'$ and $T\geq Ci=2C'i$.
If $y\in B(x,T/2)$ we have that 
$\iota$ maps $\Omega^{C'i}(X,x,x)$ into
$\Omega^{C'i + T/2}(\M,x,y)$. It follows that $\iota\circ\eta$ is a
cycle in $\Omega^{T}(\M,x,y)$. Let $[\xi]$ be the class
represented by $\iota\circ\eta$ in $\Omega (\M,x,y)$. Then,
of course, $[\xi]=I(c)$ and the lemma is proved.

\qed

We should add that I. Babenko showed in \cite{Ba1} that if 
the rational homology
 of the loop space of $X$ grows exponentially, then
$R=R_{\Omega X,0}$ where $R$ is the radius of convergence of:
\[\sum_{i\geq 2}\dim \,(\pi_{i}(X)\otimes\Q)t^{i}.\]

\begin{Corollary} Let $M$ be a closed manifold and $\M$ be its
universal covering. Suppose that there exists a 
simply connected submanifold with
boundary $X$ in $\M$ such that the map $I:H_* (\Omega (X,x,z),k_p)
\rightarrow H_* (\Omega (\M ,x,y),k_p)$ is injective (for some
$x,z\in X$, $y\in \M$ and a path from $z$ to $y$) and the
homology $H_i (\Omega X,k_p)$ grows exponentially
for some coefficient field $k_p$. Then $M$ does
not admit a Riemannian metric with vanishing topological entropy.
\label{usoenD}
\end{Corollary}

The most natural condition under which one knows that the map $I$ is
a monomorphism is when there exists a retraction from $\M$ to $X$.
Therefore we have:

\begin{Corollary} Let $M$ be a closed manifold and $\M$ be its
universal covering. Suppose that there exists a submanifold with
boundary $X$ in $\M$ such that 
$X$ is a retract of $\M$ and the 
homology $H_i (\Omega X,k_p)$ grows exponentially
for some coefficient field $k_p$. Then $M$ does
not admit a Riemannian metric with vanishing topological entropy.
\label{utilparadespues}
\end{Corollary}

\subsection{A consequence}

We now want to apply the previous results to give explicit topological
obstructions to the existence of metrics with vanishing topological
entropy. The known results up to now deal with manifolds which
have big fundamental groups (groups with exponential growth)
or manifolds with very small (finite) fundamental group (if
the loop space homology with coefficients in any field grows exponentially, 
they do not admit such a metric). We will
give an obstruction which can be applied to a great variety
of examples in the middle range (manifolds with infinite
fundamental groups with subexponential growth).

\medskip

\noindent {\bf Theorem D.} {\it Let $M$ be a closed manifold of dimension $n\geq 3$.
Suppose that $M$ can be decomposed as $X_{1}\# X_{2}$, where the order
of the fundamental group of $X_1$ is greater than 2.
If $M$ admits a $C^{\infty}$ Riemannian metric with zero
topological entropy, then $X_2$ is a homotopy sphere.}

\medskip

\begin{proof} First note that the fundamental group of $M$ must have
subexponential growth.
It is a fact from combinatorial group theory (which
follows immediately from the existence of normal forms for free
products, for instance) that if $A$ and $B$ are two finitely 
generated
groups, then the free product $A*B$ contains a free subgroup of rank
two unless $A$ is trivial or $B$ is trivial, or $A$ and $B$ are both
of order two. Since the fundamental group of a connected sum is
the free product of the fundamental groups of the summands, we
conclude that $X_2$ must be simply connected, since the 
order of the fundamental
group of $X_1$ is at least 3.

Let $\widetilde{X_1}$ be the universal covering of $X_1$. 
The universal covering $\M$ of $M=X_1 \# X_2$ is 
identified with the connected sum of $\widetilde{X_1}$ with one copy
of $X_2$ for each element of the fundamental group of
$X_1$. 
We now want to find a submanifold of $\M$ to play the role of $X$
in Corollary \ref{usoenD}.
If the fundamental group of $X_1$ is finite then we will take 
$X=\M$. Assume now that the fundamental group of $X_1$ is infinite.
Take any small $n$-dimensional disc $D$ in $\widetilde{X_1}$. Then $X=D\# 3X_2$
appears as a submanifold of $\M$. Moreover, we  claim that in
this case the map $I$ constructed as before is a monomorphism. To
prove this assume that it is not the case. Then there exists a cycle
defining a non-trivial homology class in $H_i (\Omega (X,x,z))$,
which bounds an $(i+1)$-cycle in $\Omega (\M ,x,y)$. 
The image of this
$(i+1)$-cycle is contained in some compact subset of $\M$. Performing 
connected sums at points away from this compact subset, we see that
this $(i+1)$-cycle can also be considered as a cycle in $\Omega (\M \#
3\overline{X_2},x,y)$. Then we would have that the map 
$H_* (\Omega (X,x,z))
\rightarrow H_* (\Omega (\M \# 3 \overline{X_2} ,x,y))$ is not a monomorphism.
But now we can construct a retraction $r:\Omega (\M \# 3 \overline{X_2},x,y)
\rightarrow \Omega (X,x,z)$ as follows: let $\hat{D}$ be an embedded disc 
in $\widetilde{X_1}$ containing $D$ so that $y\in \hat{D}$ and which is not
contained in the $(i+1)$-cycle mentioned before. Of course, $\hat{D} -D$
is an annulus and there is a radial retraction $\rho :\hat{D}
\rightarrow D$ sending the boundary of $\hat{D}$ to the center of $D$.
We can also assume that $\rho$ sends $y$ to $z$ and the points where
one performs the surgery with $\overline{X_2}$ to the points in $D$ 
where one performs surgery with $X_2$. The map $\rho$ then gives a
retraction $\hat{\rho} :\hat{D} \# 3 \overline{X_2} \# 3X_2
\rightarrow D\# 3X_2$ (note that the original retraction $\rho$ 
is orientation reversing outside $D$). The map $\hat{\rho}$ extends
to a retraction from $\M \# 3\overline{X_2}$ to $X$
by sending everything else to the center of $D$. But this implies that
the map  $H_*( \Omega (X,x,z))
\rightarrow H_* (\Omega (\M \# 3\overline{X_2},x,y))$ is a
monomorphism.

Up to now we have constructed in both cases (when the fundamental
group of $\M$ is finite or infinite) a submanifold $X$ of $\M$ such
that the corresponding map $I$ in the homology of the path spaces
is a monomorphism. To apply Corollary \ref{usoenD} it only remains to show
that if $X_2$ is not a homotopy sphere,
in both cases, the homology of the path spaces grow exponentially
for some coefficient field. In the second case $X$ is the connected
sum of a disc with 3 copies of $X_2$. Then the exponential growth
follows from the following lemma. In the  first case $X$ is the 
connected sum of the universal covering of $X_1$ with $k\geq 3$ copies 
of $X_2$; in this case the exponential growth follows from the
work of P. Lambrechts \cite[Theorem 3]{pascal}.

\end{proof}

\begin{Lemma} Let $M^n$ be a closed simply connected manifold which 
is not a homotopy sphere and let $N^n$ be a simply connected
compact manifold
with non-empty boundary. For a prime $p$ let $k_p$ be the field
of integers modulo $p$. 
Let $X$ be the 
compact manifold with boundary obtained
by taking the connected sum of $N$ 
with $k\geq 3$ copies of $M$. Then there exists
a prime $p$ (or $p=0$) such that the Betti numbers of
the free loop space
$\Lambda (X)$ with coefficients in $k_p$ grow exponentially.
Therefore the Betti numbers of $\Omega (X)$ also grow
exponentially.
\end{Lemma}

\begin{proof} Let us first consider the free loop space
$\Lambda (X)$. Our proof will be based
on the work of Lambrechts in \cite{pascal} and we will
will make frequent use of this reference.

Given a differential graded algebra (DGA) $(A,d)$ 
we will consider the Hochschild homology,
$HH_* (A,d)$, of $(A,d)$. We will only 
be interested in the case when the DGA is positively 
or negatively graded and connected ($A_0=R$, where
$R$ is the ground field). In this case the Hochschild 
homology is also positively or negatively graded. If 
$K$ is a simply connected $CW$-complex of finite type
and $C^* (K,k_p )$ denotes the singular cochain complex
of $K$ with coefficients in $k_p$ then 

$$H^* (\Lambda K, k_p ) \equiv HH_* (C^* (K,k_p ))$$

\noindent
(see \cite{Jones} or \cite{Halperin}). A quasi-isomorphism
(or quism) between two DGAs is a DGA-morphism which induces an
isomorphism in homology. Two DGAs are said to be
weakly equivalent if there is a chain of quasi-isomorphisms
connecting them. Since a quasi-isomorphism induces an isomorphism
between the corresponding Hochschild homologies, two weakly
equivalent DGAs have isomorphic Hochschild homologies. In summary,
to prove our lemma it is enough to obtain a DGA $(A,d)$ which is
weakly equivalent to $C^* (X, k_p )$ and such that $HH_* (A,d)$
grows exponentially (c.f. \cite[Proposition 8]{pascal}).

Since $M^n$ is not a homotopy sphere there exists a
number $i_0$, $1<i_0 <n$, such that $\pi_{i_0} (M)$ is
the first non-trivial homotopy group of $M$. Then
$\pi_{i_0} (M)=H_{i_0} (M,\Z)$. 
Actually $i_0 \leq n/2$ and so $i_0 <n-1$. 
If $Y$ is the connected sum of $k$ copies of $M$
then $H_{i_0} (Y, {\Z}) =kH_{i_0} (M,{\Z})$
(with $k\geq 3$). Therefore there exists a prime $p$ such that 
$H_{i_0} (Y, k_p )$ is a $k_p$-vector space of dimension
$\geq k$. From now on we fix the integer $i_0$ and
the prime $p$.

Pick any non-zero element in $H_{i_0} (M, k_p )$. This element gives
in a canonical way three linearly independent elements in
$H_{i_0} (Y,k_p )$ for which the cup product between any pair
of them is 0. Let $(A,d)$ be a $(i_0 -1)$-connected DGA weakly 
equivalent to the singular cochain complex of $Y$. We can
moreover assume that $d$ is decomposable
and that $A^m =0$ for all $m>n$. We can also
construct $(A,d)$ in such a way that it has three homogeneous
cycles  $x,x',y$ of degree $i_0$ (which correspond to the
cohomology classes mentioned above) such that the product of
any pair of them is 0. See for instance the argument in  
\cite[Proposition 17]{pascal}. By Poincar\'e duality there
exists  in any DGA weakly equivalent to the cochain complex
of $Y$ with $k_p$ coefficients an element of degree $n$
which is a representative of the fundamental class of $Y$ 
which can be expressed as a product of homogeneous elements
of positive degree. Let $D$ be a small open ball embedded
in $Y$.
It follows  that one can construct as
in \cite[Proposition 18]{pascal} a DGA $(\hat{A} ,\hat{d} )$
weakly equivalent to the singular cochain 
complex of $Y-D$ such that :

\vspace{.3cm}

1) $\hat{A} =A\bigoplus k_p {\bf x}$ as vector spaces, where 
$deg({\bf x})=n-1$.

2) the differential and product of $\hat{A}$ extend those of $A$.

3) ${\bf x}.{\hat{A}}^+ =0$.

4) $\hat{d} {\bf x} \in A^n$ is decomposable and represents
the fundamental class of $Y$.

\vspace{.4cm}

Now we will find a suitable DGA weakly isomorphic to
$C^* (N-D,k_p )$, where $D$ is some small open ball embedded in
$N$. Pick a 1-connected DGA $(B,d)$ weakly isomorphic to
$C^* (N,k_p )$ such that $B^m =0$ for all $m>n$. Then we consider 
the DGA $(\hat{B} ,\hat{d})$ defined as follows: as graded
vector spaces $\hat{B} =B\oplus k_p {\bf z}$ where ${\bf z}$
is a new element of degree $n-1$. The product is defined by 
extending the product of $B$ and declaring ${\bf z}.{\bf z} =0$
and ${\bf z}.B^+ =0$. It is easy to check that these choices give
a well define graded algebra.
To define the differential we set $\hat{d} =d$ on $B$ and
$\hat{d} ({\bf z})=0$. The whole idea is that ${\bf z}$ will
represent the new $(n-1)$-cohomology class represented by
$\partial D$.

Now we have to prove that $(\hat{B} ,\hat{d})$ is weakly
isomorphic to $C^* (N-D,k_p )$.

Let $\varphi : (TV,d) \rightarrow C^* (N,k_p )$ be a minimal
$TV$-model for $C^* (N,k_p )$. Let $\hat{V} =V \oplus k_p \ 
\omega \oplus W$ be a graded vector space where $\omega$ has
degree $n-1$ and the elements of $W$ are combinations of
homogeneous elements of degree $\geq n$. On $T\hat{V}$ we
define a differential $\hat{d}$ by $\hat{d}=d$ on $TV$,
$\hat{d} (\omega )=0$ and $W$ and $\hat{d}|_W$ are defined
in a minimal way to kill all cohomology in dimensions
higher than $n$. Then we extend the map $\varphi$ to
a  map

$$\hat{\varphi} :T\hat{V} \rightarrow C^*(N-D,k_p )$$

\noindent
by setting $\hat{\varphi}=i^* \varphi $ on $TV$, where $i$
is the inclusion of $N-D$ in $N$ and $\varphi (\omega )$ is
a representative of the new $(n-1)$-cohomology class 
{\it represented by} $\partial D$. It is easy to check that
$\hat{\varphi}$ is uniquely determined by this data and
it is a quasi-isomorphism. 

Since $(TV,d)$ is weakly isomorphic to $(B,d)$ there exists a
quasi-isomorphism $\psi :(TV,d)\rightarrow (B,d)$. Now define

$$ \hat{\psi} : (T\hat{V},\hat{d}) \rightarrow (\hat{B},\hat{d})$$

\noindent
by $\hat{\psi} =\psi$ on $TV$, $\hat{\psi}=0$ on $W$ and 
$\hat{\psi} (\omega )={\bf z}$. It is easy to check that
$\hat{\psi}$ is a quasi-isomorphism. Therefore $(\hat{B},
\hat{d})$ is weakly-equivalent to $C^* (N-D,k_p )$.
Moreover, $\hat{d}$ is decomposable.

Summarizing, up to now we have constructed a DGA $(\hat{A},\hat{d})$
weakly equivalent to $C^* (Y-D,k_p )$ with $\hat{d}$ decomposable
and a DGA $(\hat{B},\hat{d})$ weakly equivalent to 
$C^* (N-D,k_p )$ with $\hat{d}$ decomposable. $X$ is obtained
by joining $Y-D$ and $N-D$ through $\partial D$. A DGA 
$(F,f)$ weakly equivalent to $C^* (X,k_p )$ can then be obtained 
from $(\hat{A},\hat{d})$ and $(\hat{B},\hat{d})$ as in
\cite[Proposition 17]{pascal}:

$$F=k_p \oplus A^+ \oplus B^+ \oplus k_p {\bf w}$$

\noindent
where $f=\hat{d}$ on $A^+$ and on $B^+$ and $f({\bf w})$ is
the decomposable element in $A^n$ representing the fundamental
class of $Y$. Therefore $f$ is decomposable and  the elements
$x,x',y \in A$ satisfy all the hypotheses of 
\cite[Proposition 9]{pascal}.
Therefore the Hochschild homology of $(F,f)$ grows exponentially
and so does $H^* (\Lambda (X), k_p )$, from the discussion at
the beginning of the proof. Finally exponential growth of the
homology of $\Lambda (X)$ implies the exponential growth
of the homology of $\Omega (X)$ for any finite $CW$-complex by
an easy argument with the spectral sequence of the fibration

$$\Omega (X) \mapsto \Lambda (X) \rightarrow X.$$

\end{proof}


\section{Geometric structures and entropy}
\label{geometric-structures}

Let us begin by recalling the definition of geometric structures
in general. A geometry is a complete simply connected Riemannian
manifold $X$ such that the group of isometries acts transitively on $X$
and contains a discrete subgroup with compact quotient.
We then say that a closed manifold $M$ admits a geometric structure
modelled on $X$ if there is a Riemannian metric on $M$ such that the
Riemannian universal covering of $M$ is $X$. Maximal geometric structures in
dimension 4 have been classified by Filipkiewicz in \cite{Fi}. Wall
studies the relationship between geometric structures and complex
structures on 4-manifolds in \cite{wall}. Given a compact complex surface
we are only interested 
on whether the underlying smooth manifold admits a geometric structure,
but not on the compatibility of the geometric structure with the 
complex structure.

In this section we will describe which maximal geometric structures admit 
models with metrics with zero entropy. We refer to \cite{wall}
for a description and details of the 4-dimensional geometries.

\begin{Proposition}If $M$ admits a geometric structure modelled on one
of $S^4$, $\C P^{2}$, $S^{2}\times S^{2}$, $S^2\times\Eu^{2}$,
$\Eu^4$, $S^{3}\times \Eu^1$, {\rm Nil}$^3\times \Eu^1$ or
$\mbox{\rm Nil}^{4}$, then $M$ admits a smooth metric $g$ with
$\h =0$.
\label{geos0en}
\end{Proposition}

\begin{proof}

\begin{enumerate}
\item $S^4$, $\C P^{2}$, $S^{2}\times S^{2}$, $S^2\times\Eu^{2}$,
$\Eu^4$ and $S^{3}\times \Eu^1$:
All the Jacobi fields in these geometries grow at most linearly (in
the case of $S^4$ and $\C P^{2}$ they are actually bounded), and hence
all the Liapunov exponents of {\it every} geodesic in $M$ are zero. It
follows from Ruelle's inequality \cite{Ru} that all the measure
entropies are zero.  Hence, by the variational principle, the
topological entropy of the geodesic flow of $M$ must be zero.
\item {\rm Nil}$^3\times \Eu^1$ and $\mbox{\rm Nil}^{4}$:
The geometry Nil$^3$ can be described as $\re^{3}$ with the metric
\[ ds^{2}=dx^{2}+dy^{2}+(dz-xdy)^{2}. \]
Here, not all the Jacobi fields grow linearly, but they certainly grow
polynomially. Again this implies that all the Liapunov exponents of
{\it every} geodesic in $M$ are zero and hence the topological entropy
of the geodesic flow of $M$ must be zero.
For the case of $\mbox{\rm Nil}^{4}$ we use a result of L. Butler
\cite[Theorem 1.3]{Butler} which asserts that for any 
lattice $\Gamma$ in $\mbox{\rm Nil}^{4}$ and any left invariant metric, the 
topological entropy of the geodesic flow of $\mbox{\rm Nil}^{4}/\Gamma$
is zero. In fact, Butler's result applies to any nilpotent Lie group that
admits a normal abelian subgroup of codimension one. For general nilpotent 
groups, the result is just not true; examples have been given by Butler himself
in \cite{Butler2}.

\end{enumerate}

\end{proof}

We now show that if $M$ admits a geometric
structure modelled on one of the remaining geometries, namely
$S^2\times {\mathbb H}^{2}$, $\Eu^2\times {\mathbb H}^{2}$,
${\mathbb H}^{3}\times\Eu^1$, $\widetilde{{\rm
SL}}_{2}\times\Eu^1$, ${\rm Sol}^{4}_{0}$, ${\rm Sol}^{4}_{m,n}$
or ${\rm Sol}^{4}_{1}$, then $M$ cannot admit a metric of zero
topological entropy.  To do this, we use the next lemma, together
with the fact described in Section \ref{mecanismos}, that if $\pi_{1}(M)$ grows
exponentially, then ${\rm h}_{\rm top}(g)>0$ for any smooth metric
$g$ on $M$.

\begin{Lemma} Let $M$ be a closed orientable 4-manifold, and suppose
that $M$ admits a geometric structure modelled on one of
$S^2\times {\mathbb H}^{2}$, $\Eu^2\times {\mathbb H}^{2}$,
${\mathbb H}^{3}\times\Eu^1$, $\widetilde{{\rm
SL}}_{2}\times\Eu^1$, ${\rm Sol}^{4}_{0}$, ${\rm Sol}^{4}_{m,n}$
or ${\rm Sol}^{4}_{1}$. Then $\pi_{1}(M)$ grows exponentially.
\label{geom-exp}
\end{Lemma}

Certainly the same is true for the geometries ${\mathbb H}^2\times
{\mathbb H}^2$, ${\mathbb H}^4$ and ${\mathbb H}^2(\C)$, but here
one has a better result: the simplicial volume of any $M$ modelled
on one of these geometries will be non-zero \cite{Gromov}.

\begin{proof}
Let $(M,g)$ be a closed Riemannian manifold. Recall that $\pi_{1}(M)$
grows exponentially if and only if $\lambda(g)>0$.

If a closed manifold $M$ admits a geometric structure modelled on
$(X,G)$, then $M$ is equipped with a locally homogeneous metric
$g$. Let us consider the geometries $S^2\times {\mathbb H}^{2}$,
$\Eu^2\times {\mathbb H}^{2}$, ${\mathbb H}^{3}\times\Eu^1$. All
of them have a hyperbolic space as a factor and using this, it is
pretty straightforward to check that $\lambda(g)>0$ for all of
them.

The geometry $\widetilde{{\rm SL}}_{2}$ admits the unit sphere
bundle of a closed surface of genus $\geq 2$ as a compact
quotient. Since the fundamental group of the latter has
exponential growth, it follows that $\lambda(g)>0$ for this
geometry. Thus $\lambda(g)>0$ as well for the geometry
$\widetilde{{\rm SL}}_{2}\times\Eu^1$.

To deal with the solvable geometries, it suffices to exhibit for
each one of them a cocompact lattice $\Gamma$ such that $\Gamma$
has exponential growth.

Let $A\in SL(3,\Z)$ and let $A_{t}\in GL(3,\re)$ be a 1-parameter
subgroup with $A_{1}=A$. Let $G$ be the semidirect product
$\re^3\ltimes _{A_{t}}\re$. Let $\Gamma_{A}\subset G$ be the
cocompact lattice given by $\Z^3\ltimes _{A}\Z$ (the ascending HNN
extension of $\Z^3$ by $A$).

If we choose $A$ such that it has characteristic polynomial
$-\lambda^3+m\la^2-n\la+1$, then $G$ will be isomorphic to ${\rm
Sol}^{4}_{m,n}$. 
We may choose for example:

\[\left(\begin{array}{ccc}

1&1&0\\
m-n&m-2&1\\
m-n&m-3&1\\

\end{array}\right).\]

We claim that $\Gamma_{A}$ has exponential growth. To see
this it suffices to use a result of J. Wolf
\cite{wolf}, which asserts that a solvable (polycylic) group has
exponential growth unless it is virtually nilpotent. Wolf also
gives a criterion to check whether this holds \cite[Proposition
4.4]{wolf}. It amounts to whether $A$ has an eigenvalue with
absolute value different from one, which is certainly the case for
our choice of $A$.

To deal with ${\rm Sol}^{4}_{0}$ note that the lattice $\Gamma_A$
exists inside ${\rm Sol}^{4}_{0}$ if $A$ has eigenvalues
$\alpha>1$, $\beta$ and $\bar{\beta}$ with $\beta\neq\bar{\beta}$.
In fact with such a choice of $A$, $\Gamma_A$ is exactly the
fundamental group of the Inoue surface $S_A$. As before $\Gamma_A$
has exponential growth by Wolf's criterion.

Finally, in the case of ${\rm Sol}^{4}_{1}$, it suffices to
consider the cocompact lattice given by the ascending HNN
extension of the Heisenberg lattice $H_{Z}$ by an automorphism
given by a hyperbolic matrix in $SL(2,\Z)$. Again, by Wolf's
criterion, the lattice has exponential growth.


\end{proof}

\section{Elliptic surfaces}

Our references for the background material in this section
about elliptic surfaces are \cite{FM, wall} and references therein.

$S$ is called an elliptic surface if there is a holomorphic map
$\pi:S\to C$ whose general fibre is an elliptic curve. We will (as usual)
assume that $S$ is relatively minimal, i.e., no fibre contains an exceptional
curve of the first kind. In fact, except when $S$ is rational this is equivalent to
minimality in the usual sense.

A classification of the possible fibres of $\pi$ was given by Kodaira \cite{K1}: the cases
are labeled $I_{k}$ ($k\geq 0$), II, III, IV, I$^*$ ($k\geq 0$), II$^*$, III$^*$, IV$^*$.
Case $I_0$ means that the fibre is a smooth elliptic curve and all other types are called
{\it singular}. Also important to us is the notion of {\it multiplicity}.
Given $p\in C$, $\pi^{-1}(p)$ is a multiple fibre if there exists an integer $m>1$
such that as a divisor $\pi^*(p)=mD$. The largest such $m$ is called the multiplicity of
the fibre. Multiple fibres can only be of type $I_{k}$ for some $k$.

We will often regard $C$ as an orbifold with a $2\pi/m_i$ cone point at each point $x_i$
corresponding to a multiple fibre of multiplicity $m_i$. Then the orbifold Euler characteristic
of $C$ is:
\[\chi^{\mbox{\rm orb}}(C)=\chi(C)-\sum_{i}(1-m_i^{-1}).\]
The structure of $S$ at a smooth multiple fibre is homeomorphic to the product of a circle
with a multiple fibre in a Seifert fibration of a 3-manifold.

Let $\{t_1,\dots,t_{n}\}$ be the set of multiple points of $C$ and suppose each $t_i$ has
multiplicity $m_i$. Given a base point $t\in C$, define the {\it orbifold fundamental group}
$\pi_{1}^{\mbox{\rm orb}}(C,t)$ as follows. Let $C_0=C-\{t_1,\dots,t_n\}$. The group $\pi_1(C_0,t)$
contains the image of the free group $F$ on letters $\ga_i$, $1\leq i\leq n$, corresponding
to loops $\ga_i$ in $C_0$ enclosing $t_i$ which are null-homotopic in $C$.
Now set $\pi_{1}^{\mbox{\rm orb}}(C,t)$ to be the quotient of $\pi_1(C_0,t)$ by the smallest normal
subgroup which contains $\ga_{i}^{m_{i}}$, $1\leq i\leq n$.

We say that an orbifold is {\it good} if it has an orbifold covering whose total space
is an orbifold with no multiple points. Otherwise, an orbifold is {\it bad}.
The bad 2-dimensional orbifolds are the 2-sphere with either one multiple point
or with two multiple points with unequal multiplicities.

Of central importance to us is the fact that $\chi(S)$, the Euler number of $S$,
is always $\geq 0$ and it vanishes if and only if there are no singular fibres.

\subsection{The case of positive Euler number}

We will need the following result \cite[Theorem 2.3]{FM}.

\begin{Theorem} Let $\pi:S\to C$ be an elliptic surface. If the Euler
number of $S$ is positive, then $\pi$ induces an isomorphism
$\pi_{1}(S)\to \pi_{1}^{\mbox{\rm orb}}(C)$. \label{fungpose}
\end{Theorem}

\begin{Lemma} (a) Let $S$ be an elliptic surface with finite fundamental group.
Then, the loop space homology of the universal covering of $S$ with rational coefficients
grows exponentially.

(b) Let $S$ be a simply connected elliptic surface and let $f_{\infty}$ denote
a regular fibre. Then, the loop space homology of $S-f_{\infty}$ with rational coefficients
grows exponentially.

\label{facil}
\end{Lemma}

\begin{proof} We will use the following fact proved by J.B. Friedlander and S. Halperin \cite[Corollary 1.3]{FrH}. Let $S$ be a 1-connected finite CW-complex.
If the loop space homology of $S$ with rational coefficients
grows subexponentially, then
\[\sum_{k\geq 1}2k\,\dim\,(\pi_{2k}(S)\otimes\Q)\leq n,\]
where $n$ is the largest integer for which $H^n(S,\Q)\neq 0$.

To prove (a) observe that the Euler characteristic of $\widetilde{S}$ is given by $12d$, where $d$ is a positive integer \cite[Proposition 3.21]{FM}. Hence $b_{2}(\widetilde{S})=12d-2>2$ and by the Hurewicz theorem $\dim\,(\pi_{2}(\widetilde{S})\otimes\Q)>2$. Using the result mentioned above it follows that $S$ has the desired property.

 Let us prove item (b). Observe first that $S-f_{\infty}$ is
 simply connected. Since $H_{1}(S)=H_{3}(S)=H_{1}(S-f_{\infty})=0$,
the Mayer-Vietoris sequence gives:
\[0\to\Z^3\to H_{2}(S-f_{\infty})\oplus \Z\to H_{2}(S)\to\Z^{3}\to \Z^2\to 0.\]
It follows that $H_{2}(S-f_{\infty})$ is a free abelian group of rank
$b_{2}(S)+1$. One can deduce as above that $S-f_{\infty}$ has the desired property.

\end{proof}

We now prove:

\begin{Theorem} Let $\pi:S\to C$ be an elliptic surface. If the Euler
number of $S$ is positive, then any $C^{\infty}$ Riemannian metric
on $S$ has positive topological entropy.
\label{poseuler}
\end{Theorem}

\begin{proof} We split the proof into several cases.
\begin{enumerate}
\item $\pi_{1}^{\mbox{\rm orb}}(C)$ is finite. By Theorem \ref{fungpose},
$\pi_{1}(S)$ is finite. By Lemma \ref{facil} and Theorem C the topological entropy
of any $C^{\infty}$ metric on $S$ is positive.

\item $\pi_{1}^{\mbox{\rm orb}}(C)$ is infinite, i.e., $\chi^{\mbox{\rm orb}}(C)\leq 0$.
If $\chi^{\mbox{\rm orb}}(C)<0$, then $C$ is hyperbolic and $\pi_{1}^{\mbox{\rm orb}}(C)$
contains a free subgroup of rank two. It follows from Theorem \ref{fungpose}
that $\pi_1(S)$ grows exponentially and hence any metric on $S$ has positive
topological entropy.
If $\chi^{\mbox{\rm orb}}(C)=0$, then $C$ has a finite orbifold covering $C_0$
which has no multiple points and is a 2-torus. 
The finite cover $C_0$
induces a finite covering of the elliptic surface that we denote by $S_0$.
By Theorem 2.16 in \cite{FM}, $S_0$ is diffeomorphic to the fibre connected
sum of $\T^2\times\T^2$ with a simply connected elliptic surface $S_1$ with no multiple fibres.

We can picture $\widetilde{S_{0}}$, the universal covering of $S_0$, as $\T^{2}\times \re^2$
fibre connected sum with infinitely many copies of $S_1$ as follows.
Consider $\T^2$ (as is customary) given by the square with vertices $(0,0)$, $(1,0)$, $(1,1)$
and $(0,1)$ with the sides identified. At each fibre $\T^2\times\{(1/2+m,1/2+n)\}$
consider the fibre connected sum with $S_1$ for all integers $m$ and
$n$.
The result is  $\widetilde{S_{0}}$.
Let $X$ be $S_1$ minus an open neighborhood of a regular fibre. Let us see that
$X$ lives inside $\widetilde{S_{0}}$ as a retract. Consider a diffeomorphism of $\re^2$ that
maps the points $(1/2+m,1/2+n)$ onto the points $(0,m)$. 
This induces a diffeomorphism
of $\widetilde{S_{0}}$. Now retract
$\re^2$ over the set $[-\frac{1}{2},\frac{1}{2}]\times\re$
in the obvious way. This induces a retraction of  $\widetilde{S_{0}}$.
Now we can ``fold''
the strip onto the square $[-\frac{1}{2},\frac{1}{2}]\times
[-\frac{1}{2},\frac{1}{2}]$ in such a way that the points $(0,m)$
all go to the origin preserving orientation.
This describes  a retraction of $\re^2$ over the square 
$[-\frac{1}{2},\frac{1}{2}]\times [-\frac{1}{2},\frac{1}{2}]$
which induces a retraction of $\widetilde{S_{0}}$ onto $X$.

By Lemma \ref{facil} and Corollary \ref{utilparadespues} the topological
 entropy of any $C^{\infty}$ metric on $S_0$ (and thus of any metric on $S$) is positive.

\end{enumerate}

\end{proof}

\subsection{The case of zero Euler number}

We will need the following result \cite[Lemma 7.3, Proposition 7.4, Proposition 7.5]{FM}.

\begin{Proposition} Let $\pi:S\to C$ be an elliptic fibration with Euler number zero.
If the base orbifold is flat or hyperbolic (i.e. $\pi_{1}^{\mbox{\rm orb}}(C)$ is infinite), then
there is an exact sequence
\[0\to \Z\oplus \Z{\buildrel i_{*}\over\rightarrow} \pi_{1}(S){\buildrel\pi_{*}\over\rightarrow} \pi_{1}^{\mbox{\rm orb}}(C)\to \{1\},\]
where $i$ is the inclusion of a general fibre of $\pi$ into $S$.
In case the base orbifold $C$ is spherical and good, then there is an exact sequence
\[\Z\to \Z\oplus \Z\to \pi_{1}(S)\to \pi_{1}^{\mbox{\rm orb}}(C)\to \{1\}.\]
In case the base orbifold $C$ is bad, then $S$ is either a ruled surface over an elliptic
base ($\pi_{1}(S)\cong\Z\oplus\Z$) or $S$ is a Hopf surface with $\pi_{1}(S)\cong \Z\oplus\Z/n\Z$ for some
integer $n\geq 1$.

\end{Proposition}

Let $\kappa$ be the Kodaira dimension.

\begin{Theorem} Let $\pi:S\to C$ be an elliptic fibration with Euler number zero.
Suppose that $C$ is a good orbifold.
The following are equivalent:
\begin{enumerate}
\item $S$ admits a smooth Riemannian metric $g$ with $\h=0$;
\item $S$ admits a geometric structure modelled on 
$S^2\times \Eu^2$  ($\kappa=-\infty$, $b_1$ even)
, $\Eu^4$ ($\kappa=0$, $b_{1}$ even), $S^{3}\times \Eu^1$ ($\kappa=-\infty$, $b_1$ odd)
or {\rm Nil}$^3\times \Eu^1$ ($\kappa=0$, $b_1$ odd).
\end{enumerate}
\label{ceroeuler}
\end{Theorem}

\begin{proof} The proof is based on the following result
due to C.T.C Wall.

\begin{Theorem}[Theorem 7.4 in \cite{wall}] An elliptic surface $S$ without singular fibres has a geometric structure compatible with the complex structure of $S$ if and only if its base $C$ is a good orbifold. The type of the structure is determined as follows:

\medskip

\center \begin{tabular}{cccc}
\hline $\kappa$  &  $-\infty$  & 0 & 1\\
\hline  &   &  &    \\
$b_{1}$ even   & $S^2\times\Eu^2$  & $\Eu^4$  & $\Eu^2\times{\mathbb H}^{2}$ \\
&\\
$b_{1}$ odd   & $S^3\times\Eu^1$  & $\mbox{\rm Nil}^{3}\times\Eu^1$  &
 $\widetilde{{\rm SL}}_{2}\times\Eu^1$ \\
\hline
\end{tabular}

\end{Theorem}

\medskip

Since we are assuming that $S$ has zero Euler number and that $C$ is a
good orbifold, Wall's theorem says that $S$ admits a geometric structure as
above.
Lemma \ref{geom-exp} ensures that if $S$ admits the geometries
$\Eu^2\times \H^2$ and $\widetilde{SL}_{2}\times\Eu^1$, then $\pi_{1}(S)$ has
exponential growth and the theorem readily follows from
Proposition \ref{geos0en}.

\end{proof}

We conclude this section by noting that elliptic surfaces with
Kodaira dimension $-\infty$ and $0$ can be completely classified
\cite[Proposition 3.23]{FM} and they are usually refer to as
``not honestly" elliptic surfaces. Hence we can easily
determine the surfaces appearing in Theorem \ref{ceroeuler}:
\begin{enumerate}

\item If $S$ admits a geometric structure modelled on 
$S^2\times \Eu^2$ ($\kappa=-\infty$, $b_1$ even), then $S$
is a ruled surface of genus $1$.

\item If $S$ admits a geometric structure modelled
on $\Eu^4$ ($\kappa=0$, $b_{1}$ even), then $S$ is torus
or a hyperelliptic surface.

\item If $S$ admits a geometric structure modelled on
$S^{3}\times \Eu^1$ ($\kappa=-\infty$, $b_1$ odd), then
$S$ is a Hopf surface.

\item If $S$ admits a geometric structure modelled on
${\rm Nil}^3\times \Eu^1$ ($\kappa=0$, $b_1$ odd), then $S$
is a Kodaira surface, or a Kodaira surface modulo a finite group.

\end{enumerate}

\section{Surfaces of Kodaira dimension $-\infty$}

\subsection{Surfaces with a global spherical shell}

\begin{Proposition} A compact surface with a global spherical shell
admits no metric with zero topological entropy.
\end{Proposition}

\begin{proof} A compact surface $S$ with a global spherical shell
is diffeomorphic to a connected sum of $S^3\times S^1$ with
$b_{2}(S)$ copies of $\overline{\C P}^2$ ($b_{2}(S)>0$).

We can regard $\widetilde{S}$, the universal covering of $S$, as $S^{3}\times \re$
connected sum with infinitely many copies of $\overline{\C P^2}$ as follows.
Consider $S^1$ as $[0,1]$ with the endpoints identified and fix a point
$x\in S^3$. At each point $(x,1/2+n)\in S^3\times\re$
consider the connected sum with $\overline{\C P^2}$ for all integers $n$.
Let $X_k$ be the connected sum of $k$ copies of $\overline{\C P^2}$ with a small
open ball around a point removed. For any $k$,
$X_k$ lives inside $\widetilde{S}$ as a retract.
If $k\geq 3$, the loop space homology of $X_k$ with rational coefficients
grows exponentially and thus, by Corollary \ref{utilparadespues}, the topological entropy
of any $C^{\infty}$ metric on $S$ is positive.

\end{proof}

\subsection{Hopf surfaces}

\begin{Proposition} Any compact Hopf surface admits a smooth metric with
zero topological entropy.
\end{Proposition}

\begin{proof} Let $H$ be a finite subgroup of $U(2)$ that acts freely
on $S^3\subset\C^2$.
According to M. Kato \cite{Kato} any compact Hopf surface $S$ is diffeomorphic
to one of the following:
\begin{enumerate}
\item $S^1\times (S^3/H)$;
\item $(S^3/H)$-bundle over $S^1$ whose transition function $S^3/H\to S^3/H$
is an involution. In fact, the bundles are diffeomorphic to
$S^1\times (S^3/H)$ divided by an action of $\Z_2$ generated by
\[([t], [q])\mapsto ([t+1/2],[u(q)])\]
where $u$ is certain unitary matrix which normalizes $H$ and hence
$[u(q)]$ is well defined.
\end{enumerate}

Hence if we endow $S^1\times S^3$ with the canonical product metric
we obtain a metric with zero entropy on $S^1\times (S^3/H)$
and since $\Z_2$ acts by isometries we also obtain a metric with
zero entropy on any Hopf surface. In fact, $S$ admits a geometric
structure modelled on $S^3\times \Eu$.

\end{proof}

\subsection{Ruled surfaces}

\begin{Proposition} A ruled surface $S$ over a Riemann surface $C$ of genus $g$
 admits a smooth metric with zero topological entropy if and only if
$g=0,1$.
\end{Proposition}

\begin{proof} It follows from the homotopy exact sequence of the
fibration that if the fundamental group of $C$ grows exponentially,
then the fundamental group of $S$ also grows exponentially.
Hence if $S$ admits a smooth metric with zero topological entropy, then
the genus of $C$ must be $\leq 1$.

From the differentiable viewpoint, there are only two $S^2$-bundles
over $S^2$: the trivial one $S^2\times S^2$ and the non-trivial one which
is diffeomorphic to $\C P^2\#\overline{\C P^2}$. We explained in
\cite{PP} how to construct smooth metrics with zero entropy on these manifolds.

Similarly, from the differentiable viewpoint, there are only two $S^2$-bundles
over $\T^2$ \cite{suwa}: a trivial bundle $E_0=\T^2\times S^2$ and a non-trivial
bundle $E_1$. Clearly $E_0$ has a smooth metric with zero entropy.

The bundle $E_1$ can be described as follows (cf. \cite[p. 310]{suwa}).
Let $\Gamma$ be the group isomorphic to $\Z_2\oplus\Z_2$ generated by the
two diffeomorphisms of $S^2\times\T^2$ given by:
\[(p,[u])\mapsto (r_{z}(p), [u+1/2]),\]
\[(p,[u])\mapsto (r_{x}(p), [u+i/2]),\]
where $r_z$ is rotation by 180 degrees around the $z$-axis and
$r_x$ is rotation by 180 degrees around the $x$-axis.
Endow $S^2\times \T^2$ with the product metric. Since $\Gamma$ acts
by isometries we see that $E_1$ admits a smooth metric with zero topological
entropy. In fact $E_1$ admits a geometric structure modelled on
 $S^2\times \Eu^2$.

\end{proof}

\subsection{Inoue surfaces with $b_2=0$}

They all have fundamental group with exponential growth and hence
any metric has positive topological entropy. This can be checked directly
from the explicit presentation of the fundamental groups given in \cite{Inoue}.
Alternatively, on account of Proposition 9.1 in \cite{wall}
the Inoue surfaces admit geometric structures modelled
on solvable groups and hence we can use Lemma \ref{geom-exp}.

\section{Entropy of complex surfaces}

We are now ready to complete the discussion of the minimal entropy problem for
compact complex surfaces.
The next proposition takes care of non-minimal surfaces.

\begin{Proposition} Let $S$ be a compact complex surface. If $S$ admits
a $C^{\infty}$ metric with zero topological entropy, then $S$ must be minimal
unless $S$ is diffeomorphic to $\C P^2\#\overline{\C P^2}$.
\label{nonminimal}
\end{Proposition}

\begin{proof} If $S$ is not minimal, it is diffeomorphic to 
$S'\#\overline{\C P^2}$. Theorem D implies that the fundamental group
of $S'$ is either trivial or $\Z_2$. Moreover, the universal covering
of $S'$ must be homeomorphic to $S^4$ or $\C P^2$, otherwise the rational
loop space homology of the universal covering of $S$ has
exponential growth (cf. proof of part (a) Lemma \ref{facil} or \cite{PP}).
Since $S$ is a complex surface, the result follows.

\end{proof}

We can now directly combine Proposition \ref{nonminimal} and Theorem B with the results of Sections 6 and 7 to obtain:

\medskip

\noindent {\bf Theorem E.} {\it Let $S$ be a compact complex surface not of K\"ahler type.
 Modulo the gap in the classification of class VII surfaces we have:
The minimal entropy of $S$ is zero and the following are equivalent:
\begin{enumerate}
\item The minimal entropy problem can be solved for $S$;
\item $S$ admits a smooth metric $g$ with $\h=0$;
\item $S$ admits a geometric structure modelled on $S^{3}\times \Eu^1$ 
or {\rm Nil}$^3\times \Eu^1$;
\item $S$ has $\kappa=-\infty,0$ and is diffeomorphic to one of the following:
a Hopf surface, a Kodaira surface, or a Kodaira surface modulo a finite group.
\end{enumerate}
}

\medskip

Similarly we can combine Proposition \ref{nonminimal}, Theorem A and the results of Sections 6 and 7 to obtain:

\medskip

\noindent {\bf Theorem F.} {\it Let $S$ be a compact complex K\"ahler surface with Kodaira dimension
$\kappa\leq 1$. We have:
The minimal entropy of $S$ is zero and the following are equivalent:
\begin{enumerate}
\item The minimal entropy problem can be solved for $S$;
\item $S$ admits a smooth metric $g$ with $\h=0$;
\item $S$ admits a geometric structure modelled on $\C P^2$, $S^2\times S^2$, $S^2\times\Eu^2$
 or $\Eu^4$ or $S$ is diffeomorphic to $\C P^2\#\overline{\C P^2}$;
\item $S$ has $\kappa=-\infty,0$ and is diffeomorphic to one of the following:
$\C P^2$, a ruled surface of genus $0$ or $1$, a complex torus or
 a hyperelliptic surface.
\end{enumerate}
}

\medskip

Note that $\C P^2\#\overline{\C P^2}$ is the only non-minimal surface that admits
a metric of zero entropy. It is also the only surface that admits a metric with
zero entropy and no geometric structure.

\section{Appendix: Proof of Theorem 4.1 when $\partial X$
 is non-empty}

\begin{proof}:
Assume that the boundary of $X$ is non-empty. We have a collar
of the boundary diffeomorphic to $\partial X \times
[0,2)$. Let $Y$ be the manifold obtained by
deleting $\partial X \times [0,1)$ from $X$.
Of course, $Y$ is diffeomorphic to $X$. Now find a finite
number of convex subsets of $X$ which cover $Y$. Call
them $V_{\alpha}$, $1\leq \alpha  \leq k_0$. Let ${\bf T}$ be
a triangulation of $Y$. For a point $p\in Y$, let $F(p)$ be
the closed cell of lowest dimension containing $p$ and let
$O(p)$ be the union of all closed cells intersecting $F(p)$.
Note that $O(p)$ is a compact subset of $Y$.
Similarly, for any subset $K\subset Y$ one can define $F(Y)$
as the union of $F(p)$ for $p\in K$ and $O(K)$ as the union
of $O(p)$ for $p\in K$. 
There is a positive number $\delta$ (depending on the covering
$\{ V_{\alpha} \}$ and $g$) such that, 
after taking some barycentric
subdivisions, we can take ${\bf T}$ so that for every subset $K$
of diameter bounded by $\delta$,   $O(K)$ is
contained in one of the $V_{\alpha}$'s.

Now for $x,y\in Y$ we define the open subset $\Omega_k$ of
$\Omega (\cup V_{\alpha} ,x,y)$ as the set of all 
paths $\omega$ in $\cup V_{\alpha}$ with 
$\omega (0)=x, \omega (1)=y$ 
so that  for all $j$ between 1 and $2^k$,

$$O(\{ \omega(j-1/2^k ),\omega (j/2^k ) \} ) \ 
\cup \  \omega [j-1/2^k ,j/2^k ]$$

\noindent
is contained in one of the $V_{\alpha}$'s.
It is easy to see that $\Omega (Y,x,y)$ is
contained in the union of the $\Omega_k$'s. 

Let $B_k$ be the set of sequences $p_0 ,...,p_{2^k}$ of points 
in $Y$ such that $p_0 =x$, $p_{2^k} =y$ and for each $j$ between
1 and $2^k$ $O(\{ p_{j-1} ,p_j \} )$ is contained
in one of the $V_{\alpha}$'s. 
Let $\Omega_k^Y \subset \Omega_k$ be the set of paths in
$\Omega_k$ for which all intermediate $2^{k} -1$ points are
in $Y$. Then $B_k$ is naturally identified with a
subset of $\Omega_k^Y$ (an element of $B_k$ uniquely 
determines a broken geodesic which sends $j/2^k$ to $p_j$) and 
it is actually a deformation retract of $\Omega_k^Y$.

Given a cycle representing a homology class in
$\Omega (X,x,y)$ it can be represented by a cycle
in $\Omega (Y,x,y)$ which is therefore contained in $\Omega_k^Y$
for some $k$. Therefore we can retract it to $B_k$.
But $B_k$ is easily identified with a subset of
$Y^{2^k -1}$. Moreover, under this identification if a point
$(p_1 ,...,p_{2^k -1})\in B_k$, then the whole 
$F(p_1 )\times ...\times F(p_{2^k -1})$ is contained in
$B_k$. This implies that ${\bf T}$ induces a cell decomposition in
$B_k$. Hence the $i$-homology class can be represented by
a combination of cells of dimension $i$. A cell in $B_k$ is
a product of cells in each coordinate. The dimension of such
a cell is the sum of the dimensions of the corresponding
cells, of course. If the total dimension is $i$ then there can
be at most $i$ cells of positive dimension. 
Since $X$ is simply connected there exists a smooth map
$f:X\rightarrow X$ which is smoothly homotopic to the
identity and which sends the union of the images of all
the geodesic segments joining vertices 
in the triangulation to a point. We can moreover assume that
$x$ and $y$ are fixed by $f$. The norm of the differential
of $f$ is bounded since $X$ is compact and 
$f$ induces a map $\hat{f} : \Omega (X,x,y)\rightarrow \Omega
(X,x,y)$ which is homotopic to the identity.
Now paths belonging to an $i$-cell of $B_k$
are formed by pieces joining vertices of the triangulation and
at most $2i$ pieces in which one of the points is not
a vertex. Under the map $\hat{f}$ the former are sent to a point
and the latter to a path of length bounded in terms of the
norm of the differential of $f$ and the diameter of the 
$V_{\alpha}$'s.
Therefore, there exists a constant $C_g$
such that the image of the $i$-skeleton of
$B_k$ is sent by $\hat{f}$ to the subset of paths with
energy bounded by $C_g \ i$. The theorem follows.

\end{proof}

\end{document}